\documentclass[12pt]{article}
\usepackage{amsfonts,amssymb,amsthm}
\begin{document}
\newtheorem{thm}{Theorem}
\newtheorem{prop}{Proposition}
\newtheorem{plain}{Definition}
\newtheorem{cor}{Corollary}
\title{Classical quasi-trigonometric $r-$matrices of Cremmer-Gervais type and
their quantization}
\date{}
\author{Yermolova-Magnusson J.
\footnote{e-mail: julia@oso.chalmers.se}\\
\small{\sl Onsala Space Observatory}\\
\small{\sl Chalmers University of Technology}\\
\small{\sl G\"oteborg, Sweden}\\
Samsonov M. \footnote{e-mail: samsonov@pink.phys.spbu.ru}\\
\small{\sl St. Petersburg State University}\\
\small{\sl Institute of Physics}\\
\small{\sl St. Petersburg, Russia}\\
Stolin A. \footnote{e-mail: astolin@math.chalmers.se}\\
\small{\sl Department of Mathematics}\\
\small{\sl University of G\"oteborg}\\
\small{\sl G\"oteborg, Sweden}} \maketitle

\begin{abstract}We propose a method of quantization of certain Lie bialgebra structures on the polynomial Lie
algebras related to quasi-trigonometric solutions of the classical Yang-Baxter equation. The method is based on
so-called affinization of certain seaweed algebras and their quantum analogues.

\end{abstract}
%\maketitle
\section{Introduction}
The aim of this paper is to propose a method of quantizing certain Lie bialgebras structures on polynomial Lie
algebras. The subject was introduced by V.G. Drinfeld in the beginning of $90^{th}$, when in his paper \cite{D}
he posed the following problem: can any Lie bialgebra be quantized? Couple of years later, two mathematicians
Etingof and Kazhdan came up with the positive answer, while the problem of finding explicit quantization
formulas remained open.

First important contribution to the complete solution of the problem were presented in \cite{KM}, \cite{ESS} and
\cite{IO}, where the authors successfully quantized a list of all quasi-triangular Lie bialgebra structures on
finite dimensional simple Lie algebras also known as the Belavin-Drinfeld list.

The infinite-dimensional cases were first brought up in \cite{KST}. However, a real break-through came with
\cite{KST1} and \cite{Sam1}and was achieved due to construction of a deformed versions of Yangians
$Y(\mathfrak{sl}_{n} )$ and quantum affine algebras $U_q (\hat{\mathfrak{sl}}_{n})$ for $n=2,3$.

In the present paper we consider the case $U_q (\hat{\mathfrak{sl}}_{n})$ for $n=4$. Our solution is based on a
$q-$version of so-called seaweed algebras.

A seaweed subalgebra of $\mathfrak{sl}_{n}$ is defined as an intersection of two parabolic subalgebras one of
which containing the Borel subalgebra $B^+$ and another $B^-$. The case when both parabolic subalgebras are
maximal was treated in connection with the study of rational solutions of the classical Yang-Baxter equation.
The results of the study, as well as a complete answer to the question; when such an algebra is Frobenius, can
be found in \cite{S}. Later on in \cite{Kir} it was found out when an arbitrary seaweed algebra is Frobenius.

As a subject of this study we choose a $\mathfrak{sl}_{4}(\mathbb{C}[t])$ Lie algebra which Lie bialgebra
structure we are going to quantize. Its coalgebra is defined by a quasi-trigonometric solution of the classical
Yang-Baxter equation, this type of solutions was first introduced in \cite{KPST}. We will begin with a brief
summery of the results obtained in this paper.

Let $\mathfrak{g}$ be a simple Lie algebra and let $\Omega\in \mathfrak{g}\otimes \mathfrak{g}$ be the quadratic
Casimir element. We say that a solution $X(z,t)$ of the classical Yang-Baxter equation is {\it
quasi-trigonometric} if it is of the form:
$$X(z,t)=\frac{t\Omega}{z-t} +p(z,t),$$
where $p(z,t)$ is a polynomial with coefficients in $\mathfrak{g}\otimes \mathfrak{g}$. Any quasi-trigonometric
solution of the classical Yang-Baxter equation defines a Lie bialgebra structure on $\mathfrak{g}[t]$ and the
corresponding Lie cobracket on $\mathfrak{g}[t]$ is given by the formula
$$\{A(t)\in\mathfrak{g}[t]\}\to \{[X(z,t),A(t)\otimes 1+1\otimes A(z)]\in \mathfrak{g}[t]\otimes \mathfrak{g}[z]\}.$$
It was proved in \cite{KPST} that for $\mathfrak{g}=\mathfrak{sl}_n$ there is a one-to-one correspondence
between quasi-trigonometric $r-$matrices and Lagrangian subalgebras of $\mathfrak{g}\oplus\mathfrak{g}$
transversal to a certain Lagrangian subalgebra of $\mathfrak{g}\oplus\mathfrak{g}$ defined by a maximal
parabolic subalgebra of $\mathfrak{g}$. Here we define Lagrangian with respect to the following symmetric
non-degenerate invariant bilinear form on $\mathfrak{g}\oplus\mathfrak{g}$:
$$Q((a,b),(c,d))=K(a,c)-K(b,d),$$
where $K$ is the Killing form.

We are going to quantize a quasi-trigonometric $r-$matrix with coefficients in
$\mathfrak{sl}_4\otimes\mathfrak{sl}_4$ defined by the following Lagrangian subalgebra $W$ of
$\mathfrak{sl}_4\oplus\mathfrak{sl}_4$:
$$W=\{(X,Y)\in \mathfrak{sl}_4\oplus\mathfrak{sl}_4:\ X=Ad(T)Y\}.$$
Here
$$
T= \left( {\begin{array}{cccc}
0 & 1 & 0 & 0\\[2ex]
1 & 0 & 1 & 0\\[2ex]
0 & 0 & 1 & 0\\[2ex]
0 & 1 & 0 & 1\\[2ex]
\end{array}
} \right).
$$
\section{Preliminaries}
We denote by $\hat{\mathfrak{sl}}_{n}$, $n>2$, the following affine Lie algebra with the set of simple affine
roots $\{\alpha_{0}, \alpha_{1},\ldots,\alpha_{n-1}\},$ defined by the generators
$\{e_{\pm\alpha_{i}},h_{\alpha_{i}}\}_{0\le i\le n-1}$ and the relations
\begin{equation}
{
\begin{array}{ccc}
[h_{\alpha_{i}},e_{\pm\alpha_{j}}]=\pm (\alpha_{i},\alpha_{j})e_{\pm\alpha_{j}}, &&
[e_{\alpha_{i}},e_{-\alpha_{j}}]=\delta_{ij}\displaystyle h_{\alpha_{i}};
\end{array}
}
\end{equation}
\begin{eqnarray}
e_{\pm\alpha_{i}}^{2}e_{\pm\alpha_{j}}-2\mathop{}
e_{\pm\alpha_{i}}e_{\pm\alpha_{j}}e_{\pm\alpha_{i}}+e_{\pm\alpha_{j}}e_{\pm\alpha_{i}}^{2}=0,&& \mbox{ if
}|i-j|=1;
\\
e_{\pm\alpha_{i}}e_{\pm\alpha_{j}}-e_{\pm\alpha_{j}}e_{\pm\alpha_{i}}=0,&& \mbox{ if }|i-j|\ne 1.
%\label{rel4}
\end{eqnarray}
Here $\{(\alpha_i,\alpha_j)\}$ is the affine Cartan matrix of the root system $A_{n-1}^{(1)}$ that is:
$\{(\alpha_i,\alpha_j)\}=2$ if $i=j$, $\{(\alpha_i,\alpha_j)\}=-1$ if $|i-j|=1,n-1$ and $0$ otherwise. The root
$\alpha_{0}+\alpha_{1}+\ldots+\alpha_{n-1}$ is called imaginary and is denoted by $\delta$, the corresponding
Cartan element $h_{\delta}=\sum_i h_{\alpha_i}$ is central in $\hat{\mathfrak{sl}}_{n}$. It is customary to
exclude $\alpha_0$ from the defining relations substituting it by $\delta-\alpha_1-\ldots-\alpha_{n-1}$.
\\

The universal enveloping algebra $U( \hat{\mathfrak{sl}}_{n})$ has the following quantum analogue denoted by
$U_q ( \hat{\mathfrak{sl}}_{n})$, which is a Hopf algebra. Formally it is defined by the same set of generators
as $U( \hat{\mathfrak{sl}}_{n})$, however the defining relations are being $q-$deformed.
\begin{equation}
{
\begin{array}{ccc}
[h_{\alpha_{i}},e_{\pm\alpha_{j}}]=\pm (\alpha_{i},\alpha_{j})e_{\pm\alpha_{j}}, &&
[e_{\alpha_{i}},e_{-\alpha_{j}}]=\delta_{ij}\displaystyle
\frac{q^{h_{\alpha_{i}}}-q^{-h_{\alpha_{i}}}}{q-q^{-1}};
\end{array}
} \label{rea1}
\end{equation}
\begin{equation}
e_{\pm\alpha_{i}}^{2}e_{\pm\alpha_{j}}-(q+q^{-1})\mathop{} e_{\pm\alpha_{i}}
e_{\pm\alpha_{j}}e_{\pm\alpha_{i}}+e_{\pm\alpha_{j}}e_{\pm\alpha_{i}}^{2}=0, \mbox{ if }|i-j|=1;
\end{equation}
\begin{equation}
e_{\pm\alpha_{i}}e_{\pm\alpha_{j}}-e_{\pm\alpha_{j}}e_{\pm\alpha_{i}}=0, \mbox{ if }|i-j|\ne 1.
\end{equation}
A Hopf algebra structure is defined uniquely by the following values of the coproduct on the Chevalley
generators
\begin{equation}
\Delta(h_{\alpha_{i}})=h_{\alpha_{i}}\otimes 1+1\otimes h_{\alpha_{i}}; \label{rel6}
\end{equation}
\begin{equation}
\Delta(e_{\alpha_{i}})=q^{-h_{\alpha_{i}}}\otimes e_{\alpha_{i}}+ e_{\alpha_{i}}\otimes 1,\hphantom{aa}
\Delta(e_{-\alpha_{i}})=e_{-\alpha_{i}}\otimes q^{h_{\alpha_{i}}}+1\otimes e_{-\alpha_{i}}; \label{rel7}
\end{equation}
\begin{equation}
\begin{array}{ccc}
S(h_{\alpha_{i}})=-h_{\alpha_{i}},& S(e_{\alpha_{i}})=-q^{h_{\alpha_{i}}} e_{\alpha_{i}},&
S(e_{-\alpha_{i}})=-e_{-\alpha_{i}}q^{-h_{\alpha_{i}}};
\end{array}
\end{equation}
\begin{equation}
\begin{array}{ccc}
\varepsilon(h_{\alpha_{i}})=0, &\varepsilon(e_{\alpha_{i}})=0, & \varepsilon(e_{-\alpha_{i}})=0.
\end{array}
\end{equation}
Following \cite{KT, KT2}, we fix a normal ordering
$$
\begin{array}{l}
\alpha_{1}\prec\alpha_{1}+\alpha_{2}\prec\cdots\prec\alpha_{1}+\alpha_{2}+\ldots+
\alpha_{n-1}\prec\alpha_{2}\prec\cdots\prec\alpha_{n-1}\prec\\
\prec\ldots\prec\delta\prec\ldots\prec\delta-\alpha_{n-1}\prec\delta-\alpha_{n-1}-\alpha_{n-2}
\prec\delta-\alpha_{1}-\ldots\prec\alpha_{n-1}.
\end{array}
$$
Let $\alpha$, $\beta$ and $\gamma$ be pairwise non-collinear roots taken so that there are no other roots
$\alpha^{\prime}$ and $\beta^{\prime}$ with the property $\gamma=\alpha^{\prime}+\beta^{\prime}.$ Then, if
$e_{\pm\alpha}$ and $e_{\pm\beta}$ have already been constructed, one sets
$$
\begin{array}{ll}
e_{\gamma}=e_{\alpha}e_{\beta}-q^{(\alpha,\beta)}e_{\beta}e_{\alpha},&
e_{-\gamma}=e_{-\beta}e_{-\alpha}-q^{-(\beta,\alpha)}e_{-\alpha}e_{-\beta}.
\end{array}
$$
Consider a seaweed Lie algebra $SW_{n}:=P_{1}^{+}\cap P_{n-1}^{-}$ which is defined as the intersection of the
maximal parabolic subalgebras $P_{i}^{\pm}$, obtained by omitting the generators corresponding to the roots
containing $\pm\alpha_{i}$. Classically, we have an embedding of the restricted seaweed Lie algebra $\iota:
SW_{n}^{\prime}\hookrightarrow\hat{\mathfrak{sl}}_{n-1}$, where
$$
SW_{n}^{\prime}:=(SW_{n}\backslash\{h_{\alpha_{1}},h_{\alpha_{n-1}}
\})\bigcup\{p=\sum_{k=1}^{n-1}\frac{(n-2k)}{n}h_{\alpha_{k}}\}.
$$
We show that this embedding can be \textbf{\textit{quantized}} in the sense that it can be lifted to an
embedding:
\begin{equation}
\iota_{q}:U_{q}(SW_{n}^{\prime})\hookrightarrow U_{q}^{{\cal K}_{n-1}}(\hat{\mathfrak{sl}}_{n-1}) \label{isom}
\end{equation}
Here the coalgebra structure on $U_{q}(SW_{n}^{\prime})$ is obtained from the standard (see above) coalgebra
structure on $U_q(\mathfrak{sl}_{n})$ by twisting it by a certain element ${\cal K}_{n}\in
U_{q}(\mathfrak{sl}_{n})^{\otimes 2}$, which will be described below.

The element ${\cal K}_{n}$ will be chosen as follows: ${\cal K}_{n}=q^{r_{0}}$, where $r_{0}$ is the Cartan part
of the so-called Cremmer-Gervais $r-$matrix. The Cremmer-Gervais $r-$matrices correspond to the
\textit{shift-by-one} Belavin-Drinfeld triples: $\tau:\Gamma_{1}\rightarrow~\Gamma_{2}$, where
$\Gamma_{i}=\{\alpha_{i},\ldots,\alpha_{n+i-2}\}, i=1,2$. To explain this in more details let us recall
classification of quasi-triangular $r-$matrices for a simple Lie algebra $\mathfrak{g}$. First of all the
quasi-triangular $r-$matrices are solutions of the system
\begin{eqnarray}
r^{12}+r^{21}&=&\Omega\label{dr1}\\[2ex]
[r^{12},r^{13}]+[r^{12},r^{23}]+[r^{13},r^{23}]&=&0 \label{dr2}
\end{eqnarray}
where $\Omega$ is the quadratic Casimir element in $\mathfrak{g}\otimes\mathfrak{g}$.

Belavin and Drinfeld proved that any solution of this system is defined by a triple $(\Gamma_{1},\Gamma_{2},
\tau)$, where $\Gamma_{1},\Gamma_{2}$ are subdiagrams of the Dynkin diagram of $\mathfrak{g}$ and $\tau$ is an
isometry between these two subdiagrams. Further, each $\Gamma_i$ defines a reductive subalgebra of
$\mathfrak{g}$, and $\tau$ extends to an isometry (with respect to the corresponding restrictions of the Killing
form) between the obtained two reductive subalgebras of $\mathfrak{g}$. The following property of $\tau$ should
be satisfied: $\tau^{k}(\alpha)\not\in\Gamma_{1}$ for any $\alpha\in\Gamma_{1}$ and some $k$. Let $\Omega_0$ be
the Cartan part of $\Omega$. Then one can construct a quasi-triangular $r-$matrix according to the following
\begin{thm}[Belavin-Drinfeld]
Let $r_{0}\in\mathfrak{h}\otimes\mathfrak{h}$ satisfies the systems
\begin{eqnarray}
r_{0}^{12}+r_{0}^{21}&=&\Omega_{0},\label{bd1}\\[2ex]
(\tau(\alpha)\otimes 1+1\otimes\alpha)(r_{0})&=&0\label{bd2}
\end{eqnarray}
for any $\alpha\in\Gamma_{1}$. Then the tensor
$$
r=r_{0}+\sum_{\alpha>0}X_{-\alpha}\otimes X_{\alpha}+ \sum_{\alpha,\beta>0,\alpha\succ\beta}X_{-\alpha} \wedge
X_{\beta}
$$
satisfies (\ref{dr1}),(\ref{dr2}). Moreover, any solution of (\ref{dr1}),(\ref{dr2}) is of the above form, for a
suitable triangular decomposition of $\mathfrak{g}$ and suitable choice of a basis $\{X_{\alpha}\}$.
\end{thm}
\section{Quantization of $U_{q}(SW_{5}^{\prime})$ and its affinization}
We define  $SW_{5}\subset \mathfrak{sl}_{5}$ is the following seaweed Lie subalgebra:
$$
\left( {
\begin{array}{ccccccccc}
\ast && 0 && 0 && 0 && 0\\
\ast && \ast && \ast && \ast && \ast\\
\ast && \ast && \ast &&  \ast && \ast\\
\ast && \ast && \ast && \ast && \ast\\
0 && 0 && 0 && 0 && \ast\\
\end{array}
} \right).
$$
Consider the restricted seaweed algebra
$$
SW^{\prime}_{5}:=(SW_{5}\backslash\{h_{\alpha_{1}}, h_{\alpha_{4}}\})\bigcup\{p_{5}:= \frac 35
h_{\alpha_{1}}+\frac 15 h_{\alpha_{2}}- \frac 15 h_{\alpha_{3}}-\frac 35 h_{\alpha_{4}}\}.
$$
\begin{prop}
There exists an embedding: $\iota^{(5)}: SW^{\prime}_{5}\hookrightarrow\hat{\mathfrak{sl}}_{4}$.
\end{prop}
\begin{proof}
We define $\iota^{(5)}$ by its values on the generators of $SW_{5}^{\prime}$ letting (to avoid ambiguities, we
mark the generators of the affine algebra $\hat{\mathfrak{sl}}_{4}$ by hat, for instance: $\hat{e}_{\alpha}$).
$$
\begin{array}{lcl}
\iota^{(5)}(p_{5})=\frac 34 h_{\alpha_{1}}+\frac 24 h_{\alpha_{2}}+ \frac 14 h_{\alpha_{3}},&&
\iota^{(5)}(e_{\alpha_{4}})=\hat{e}_{\delta-\alpha_{1}-\alpha_{2}-\alpha_{3}},\\[2ex]
\iota^{(5)}(e_{-\alpha_{1}})=\hat{e}_{-\alpha_{1}},&&
\iota^{(5)}(e_{\alpha_{2}})=\hat{e}_{\alpha_{2}},\\[2ex]
\iota^{(5)}(e_{\alpha_{3}})=\hat{e}_{\alpha_{3}}, &&
\iota^{(5)}(h_{\alpha_{2}})=\hat{h}_{\alpha_{2}},\\[2ex]
\iota^{(5)}(h_{\alpha_{3}})=\hat{h}_{\alpha_{3}}.&&
\end{array}
$$
and fixing the normal ordering in $\hat{\mathfrak{sl}}_{4}$ so that
$\alpha_{1}\prec\alpha_{2}\prec\alpha_{3}\prec \delta-\alpha_{1}-\alpha_{2}-\alpha_{3}$. The proof will be
finished once we will verify that the defining relations of $SW_{5}^{\prime}$ are being preserved by
$\iota^{(5)}$. This is straightforward.
\end{proof}
We define $U_q (SW'_5)$ as a subalgebra of $U_q (\mathfrak{sl}_5)$ generated by $h_{\alpha_2},$ $ h_{\alpha_3},
p_5,$ $ e_{-\alpha_1},$ $ e_{\alpha_4}, $ $e_{\pm\alpha_2}, $ $e_{\pm\alpha_3}.$ Unfortunately, the subalgebra
defined as $U_q (SW'_{5})$ is not a Hopf subalgebra of $U_{q}(\mathfrak{sl}_{5})$ with respect to the standard
comultiplication in $U_{q}(\mathfrak{sl}_{5})$. However, twisting the standard comultiplication by  ${\cal
K}_{5}=q^{r_{0}(5)}$, where $r_0 (5)$ is   the Cartan part of the Cremmer-Gervais $r-$matrix for
$\mathfrak{sl}_5$, we obtain the following result:
\begin{prop}
Let ${\cal K}_{5}=q^{r_{0}(5)}$ where
$$
\begin{array}{l}
r_{0}(5)=\\[2ex]
h_{\alpha_{1}}\otimes (\frac 25 h_{\alpha_{1}}+ \frac 35 h_{\alpha_{2}}+ \frac 35 h_{\alpha_{3}}+\frac 25
h_{\alpha_{4}})+h_{\alpha_{2}}\otimes( \frac 35 h_{\alpha_{2}}+\frac 45 h_{\alpha_{3}}+\frac 35 h_{\alpha_{4}})+
\\[2ex]
+h_{\alpha_{3}}\otimes (-\frac 15 h_{\alpha_{1}}+\frac 35 h_{\alpha_{3}} +\frac 35 h_{\alpha_{4}})
+h_{\alpha_{4}}\otimes(-\frac 15 h_{\alpha_{1}}-\frac 15 h_{\alpha_{2}} +\frac 25 h_{\alpha_{4}})
\end{array}
$$
is the Cartan part of the Cremmer-Gervais $r-$matrix. Then the Hopf algebra $U_{q}^{{\cal
K}_{5}}(\mathfrak{sl}_{5})$ induces a coalgebraic structure on the quantum seaweed algebra
$U(SW_{5}^{\prime})\subset U_{q}^{{\cal K}_{5}}(\mathfrak{sl}_{5})$.
\end{prop}
\begin{proof}
Let us consider  $U_{q}^{{\cal K}_{5}}(\mathfrak{sl}_{5})$ and check that the elements corresponding to the
generators of $U_{q}(SW_{5}^{\prime})$ define a Hopf subalgebra in $U_{q}^{{\cal K}_{5}}(\mathfrak{sl}_{5})$. In
order to prove this, we have to show that the following relations hold for $1\le i\le 4$:
\begin{eqnarray}
(\alpha_{i}\otimes{\rm id})(r_{0}(5))&\in &U_{q}(SW_{5}^{\prime})
\label{est1}\\
-h_{\alpha_{i}}+({\rm id}\otimes\alpha_{i})(r_{0}(5)))
&\in &U_{q}(SW_{5}^{\prime})\label{est2}\\
({\rm id} \otimes\alpha_{i})(r_{0}(5))&\in & U_{q}(SW_{5}^{\prime})\label{est3}\\
h_{\alpha_{i}}+(\alpha_{i}\otimes{\rm id})(r_{0}(5)) &\in & U_{q}(SW_{5}^{\prime}).\label{est4}
\end{eqnarray}
Taking into account that $h_{\alpha_{2}},h_{\alpha_{3}}\in U_{q}(SW^{\prime}_{5})$ and formulas
(\ref{bd1})-(\ref{bd2}), we see that using the explicit formula for $r_{0}(5)$ it is sufficient to prove that
\begin{equation}
({\rm id}\otimes\alpha_{1})(r_{0}(5))=\frac{1}{5} h_{\alpha_{1}}-\frac{3}{5}h_{\alpha_{2}}-\frac 25
h_{\alpha_{3}} -\frac{1}{5}h_{\alpha_{4}}\in U_{q}(SW_{5}^{\prime}). \label{bd5}
\end{equation}
It can be checked directly and this finishes the proof.
\end{proof}
In what follows we need to introduce an element ${\hat{{\cal K}}_{4}}\in U_{q}(\hat{\mathfrak{sl}}_{4})^{\otimes
2}$, which is a four-dimensional analogue of the element
${\cal K}_{5}\in U_{q}({\mathfrak{sl}}_{5})^{\otimes 2}$:\\
$\hat{{\cal K}}_{4}=q^{\hat{r}_0(4)}$, where
$$\begin{array}{l}
\hat{r}_{0}(4)=\\[2ex]
\hat{h}_{\alpha_{1}}\otimes (\frac 38\hat{h}_{\alpha_{1}}+ \frac 12 \hat{h}_{\alpha_{2}}+ \frac 38
\hat{h}_{\alpha_{3}})+\hat{h}_{\alpha_{2}}\otimes( \frac 12\hat{h}_{\alpha_{2}}+\frac 12 \hat{h}_{\alpha_{3}})
+\hat{h}_{\alpha_{3}}\otimes (-\frac 18 \hat{h}_{\alpha_{1}}+\frac 38 \hat{h}_{\alpha_{3}} ).
\end{array}
$$
\begin{thm}
There is an embedding of Hopf algebras:
$$
\iota_{q}^{(5)}: U_{q}(SW^{\prime}_{5}) \hookrightarrow U_{q}^{\hat{{\cal K}}_{4}}(\hat{\mathfrak{sl}}_{4}),
$$
where $\iota_{q}^{(5)}$ is defined by the same formulas on the quantum generators as $\iota^{(5)}$ was defined
on the classical generators in Proposition 1.
\end{thm}
\begin{proof}
By the properties (\ref{dr1})-(\ref{dr2}) it is sufficient to verify that
\begin{eqnarray}
\iota_{q}^{(5)}(({\rm id}\otimes\alpha_{1})(r_{0}(5)))
&=&({\rm id}\otimes\alpha_{1})(\hat{r}_{0}(4))\label{sv3}\\
\iota_{q}^{(5)}((\alpha_{4}\otimes {\rm id}) (r_{0}(5)))&=&-((\alpha_{1}+\alpha_{2}+\alpha_{3}) \otimes{\rm
id})(\hat{r}_{0}(4))\label{sv4}
\end{eqnarray}
(\ref{sv3}) is equivalent to the following equality
$$
\iota_{q}^{(5)}\left(\frac{1}{5} h_{\alpha_{1}}-\frac{3}{5} h_{\alpha_{2}}-\frac{2}{5} h_{\alpha_{3}}-\frac 15
h_{\alpha_{4}}\right)= \frac{1}{4}\hat{h}_{\alpha_{1}}-\frac 12
\hat{h}_{\alpha_{2}}-\frac{1}{4}\hat{h}_{\alpha_{3}}
$$
while (\ref{sv4}) can be rewritten as follows:
$$
\iota_{q}^{(5)}\left(-\frac{1}{5}h_{\alpha_{1}}- \frac 25 h_{\alpha_{2}}-\frac 23 h_{\alpha_{3}}+\frac 15
h_{\alpha_{4}}\right) =-\frac 14 \hat{h}_{\alpha_{1}}-\frac 12 \hat{h}_{\alpha_{2}}-\frac 34
\hat{h}_{\alpha_{3}}.
$$
\end{proof}
\section{Quantization of quasi-trigonometric $r-$ matrices}
Following \cite{ESS, IO}, we can construct the Cremmer-Gervais twist in $U_{q}(\mathfrak{sl}_{5})$ as the
following product:
$$
{\cal F}_{CG_{5}}={\cal F}^{(3)}{\cal F}^{(2)}{\cal F}^{(1)}\cdot {\cal K}_{5}
$$
where
$$
{\cal
F}^{(k)}=\prod_{\alpha_{1}\preceq\beta\preceq{\alpha_{3}}}\exp_{q^{2}}((q^{-1}-q)e_{\tau^{k}(\beta)}^{\prime}
\otimes e_{-\beta}^{\prime}),
$$
$$
\exp_{q}(x):=\sum_{n\ge 0}\frac{x^{n}}{(n)_{q}!}, \hphantom{aa} (n)_{q}!\equiv (1)_{q}(2)_{q}\ldots
(n)_{q},\hphantom{aa} (k)_{q}\equiv (1-q^{k})/(1-q)
$$
and
$$
\begin{array}{lcl}
e_{\alpha_{k}}^{\prime}=q^{(\alpha_{k}\otimes{\rm id}) (r_{0}(5))}e_{\alpha_{k}},&&
e_{-\alpha_{k}}^{\prime}=q^{-({\rm id}\otimes\alpha_{k}) (r_{(0)}(5))}e_{-\alpha_{k}}.
\end{array}
$$
The elements $e_{\pm\alpha_{k}}^{\prime}$ have the following coproducts after twisting by ${\cal K}_{5}$:
$$
\begin{array}{lcl}
{\cal K}_{5} \Delta(e_{\alpha_{k}}^{\prime}) {\cal K}_{5}^{-1}&=& e_{\alpha_{k}}^{\prime} \otimes
q^{2(\alpha_{k}\otimes{\rm id})(r_{0}(5))}+
1\otimes e_{\alpha_{k}}^{\prime}\\[2ex]
{\cal K}_{5}\Delta(e_{-\alpha_{k}}^{\prime}) {\cal K}_{5}^{-1}&=&q^{-2({\rm id}\otimes\alpha_{k})(r_{0}(5))}
\otimes e_{-\alpha_{k}}^{\prime}+e_{-\alpha_{k}}^{\prime} \otimes 1.
\end{array}
$$
\begin{cor}
${\cal F}^{\prime}_{CG_{5}}={\cal F}_{CG_{5}}{\cal K}_{5}^{-1}$ restricts to a twist of $U_{q}(SW_{5}^{\prime})$
leading to  an affine twist in $U_{q}(\hat{\mathfrak{sl}}_{4})$ via embedding $\iota_{q}^{(5)}$:
\begin{equation}
\hat{{\cal F}}_{CG_{4}}=(\iota_{q}^{(5)}\otimes\iota_{q}^{(5)}) ({\cal F}_{CG_{5}}\cdot{\cal K}_{5}^{-1})
\cdot\hat{{\cal K}}_{4}. \label{affinization}
\end{equation}
\end{cor}
\begin{proof}
Since the element ${\cal K}_{5}$ is a twist for $U_{q}(\mathfrak{sl}_{5})$ we see that:
$$
{\cal K}_{5}^{12}(\Delta\otimes{\rm id})({\cal K}_{5}) ={\cal K}_{5}^{23}({\rm id}\otimes\Delta)({\cal K}_{5}),
$$
then, by the chain property of the twists, ${\cal F}^{\prime}_{CG_{5}}$ defines a twist for $U_{q}^{{\cal
K}_{5}}(\mathfrak{sl}_{5})$, and it can be restricted to a twist for $U_{q}(SW'_{5})$ leading to
(\ref{affinization}) if we use the same argument relating $U_{q}(\hat{\mathfrak{sl}}_{4})$ and
$U_{q}^{\hat{{\cal K}}_{4}}(\hat{\mathfrak{sl}}_{4})$ via twisting by $\hat{{\cal K}}_{4}$.
\end{proof}
\section{Discussion}
It follows from Proposition 1 that in the present paper we have quantized the following quasi-trigonometric
$r-$matrix in $\mathfrak{sl}_4$:
$$
\begin{array}{l}
X(z,t)=\displaystyle\frac{t\Omega}{z-t}+r_{DJ}+\\[2ex]
SK((e_{23}+e_{34})\otimes e_{21}+e_{24}\otimes e_{31}+e_{34}\otimes e_{32}-r_0(4)+\\[2ex]
(z-t)(e_{21}\otimes e_{41}+e_{41}\otimes e_{21}+e_{31}\otimes e_{31})+\\[2ex]
ze_{31}\otimes e_{42}- te_{42}\otimes e_{31} +ze_{41}\otimes e_{32}-te_{32}\otimes e_{41}.
\end{array}
$$
Here $e_{ij}$ are matrix units, $SK(a\otimes b)=a\otimes b-b\otimes a$ and $r_{DJ}$ is the Drinfeld-Jimbo
$r-$matrix for $\mathfrak{sl}_4$ In terms of Theorem 1 it is defined by the following data:
$$\Gamma_1=\Gamma_2=\emptyset,\ \ r_0=\displaystyle\frac{\Omega_0}{2}$$
The corresponding quantum algebra can be obtained by twisting the standard quantum affine algebra
$U_q(\hat{\mathfrak{sl}}_4)$ by the element $\hat{\cal{F}}_{CG_4}$. Then the quantum $R-$matrix
$\hat{\cal{F}}_{CG_4}^{21}R_{aff}\hat{\cal{F}}_{CG_4}^{-1}$ quantizes $X(z,t)$. The universal quantum $R-$matrix
$R_{aff}$ was constructed in \cite{KT2}.

Using methods of \cite{KPST} one can prove that this $r-$matrix is gauge equivalent to that of defined by the
Lagrangian subalgebra $W\subset \mathfrak{sl}_4\oplus \mathfrak{sl}_4$ from Introduction. Extra terms of
$X(z,t)$,
$$X(z,t)-\displaystyle\frac{t\Omega}{z-t}-r_{DJ},$$
define a classical twist for the polynomial Lie bialgebra $\mathfrak{sl}_4 [t]$.

As a conclusion, we've to mention that the results presented in this paper support the conjecture made in
\cite{KPST}, namely:

{\it every classical twist gives rise to a quantum twist.}

\end{document}